\documentclass{amsart}
\usepackage{amsthm,amsfonts,amsmath}
\usepackage[numbers,sort&compress]{natbib}

\theoremstyle{plain}
\newtheorem{theorem}{Theorem}
\newcommand{\thmref}[1]{Theorem~\ref{thm:#1}}
\newcommand{\secref}[1]{Section~\ref{sec:#1}}
\newcommand{\seclabel}[1]{\label{sec:#1}}
\newcommand{\thmlabel}[1]{\label{thm:#1}}
\newcommand{\etal}{ \emph{et al.}~}
\newcommand{\floor}[1]{\ensuremath{\protect\lfloor#1\rfloor}}
\newcommand{\tpw}[1]{\ensuremath{\textup{\textsf{tpw}}(#1)}}
\newcommand{\tw}[1]{\ensuremath{\textup{\textsf{tw}}(#1)}}

\renewcommand{\baselinestretch}{1.175}
\begin{document}

\title{Vertex Partitions of Chordal Graphs}

\author{David R. Wood}
\address{School of Computer Science, Carleton University, Ottawa, Canada.
Department of Applied Mathematics, Charles University, Prague, Czech Republic.}
\email{davidw@scs.carleton.ca}
\thanks{Research supported by NSERC and COMBSTRU}
\date{\today}


\begin{abstract} 
A \emph{$k$-tree} is a chordal graph with no $(k+2)$-clique. An \emph{$\ell$-tree-partition} of a graph $G$ is a vertex partition of $G$ into `bags', such that contracting each bag to a single vertex gives an $\ell$-tree (after deleting loops and replacing parallel edges by a single edge). We prove that for all $k\geq\ell\geq0$, every $k$-tree has an $\ell$-tree-partition in which every bag induces a connected $\floor{k/(\ell+1)}$-tree. An analogous result is proved for oriented $k$-trees.
\end{abstract}


\maketitle

\section{Introduction}
\seclabel{Introduction}


Let $G$ be an (undirected, simple, finite) graph with vertex set $V(G)$ and edge set $E(G)$. The neighbourhood of a vertex $v$ of $G$ is denoted by $N(v)=\{w\in V(G):vw\in E(G)\}$. A \emph{chord} of a cycle $C$ is an edge not in $C$ whose endpoints are both in $C$. $G$ is \emph{chordal} if every cycle on at least four vertices has a chord. A $k$-clique ($k\geq0$) is a set of $k$ pairwise adjacent vertices. A \emph{$k$-tree} is a chordal graph with no $(k+2)$-clique. The \emph{tree-width} of $G$, denoted by $\tw{G}$, is the minimum $k$ such that $G$ is a subgraph of a $k$-tree. It is well known that $G$ is a $k$-tree if and only if $V(G)=\emptyset$, or $G$ has a vertex $v$ such that $G\setminus v$ is a $k$-tree, and $N(v)$ is a $k'$-clique for some $k'\leq k$.


Let $G$ and $H$ be graphs. The elements of $V(H)$ are called \emph{nodes}. 
Let $\{H_x\subseteq V(G):x\in V(H)\}$ be a  set of subsets of $V(G)$ indexed by the nodes of $H$. Each set $H_x$ is called a \emph{bag}. The pair $(H,\{H_x\subseteq V(G):x\in V(H)\})$ is an \emph{$H$-partition} of $G$ if:
\begin{itemize}
\item $\forall$ vertices $v$ of $G$, $\exists$ node $x$ of $H$ with $v\in H_x$, and
\item $\forall$ distinct nodes $x$ and $y$ of $H$, $H_x\cap H_y=\emptyset$, and
\item $\forall$ edge $vw$ of $G$,  either 
\begin{itemize}
\item $\exists$ node $x$ of $H$ with $v\in H_x$ and $w\in H_x$, or
\item $\exists$ edge $xy$ of $H$ with $v\in H_x$ and $w\in H_y$.
\end{itemize}
\end{itemize}

For brevity we say $H$ is a partition of $G$. A \emph{$k$-tree-partition} is an $H$-partition for some $k$-tree $H$. A \emph{tree-partition} is a $1$-tree-partition. Tree-partitions were independently introduced by \citet{Seese85} and \citet{Halin91}, and have since been investigated by a number of authors \cite{BodEng-JAlg97, Bodlaender-DMTCS99, Halin91, Seese85, DO-JGT95, DO-DM96}. The main property of tree-partitions that has been studied is the maximum cardinality of a bag, called the \emph{width} of the tree-partition. The minimum width over all tree-partitions of a graph $G$ is the \emph{tree-partition-width}\footnote{Tree-partition-width has also been called \emph{strong tree-width} \cite{Seese85,BodEng-JAlg97}.} of $G$, denoted by $\tpw{G}$. A graph with bounded degree has bounded tree-partition-width if and only if it has bounded tree-width \cite{DO-DM96}. In particular, for every graph $G$, \citet{Seese85} proved that $\tw{G}\leq2\,\tpw{G}-1$, and \citet{DO-JGT95} proved that $\tpw{G}\leq24\,\tw{G}\max\{\Delta(G),1\}$, where $\Delta(G)$ is the maximum degree of $G$. See \citep{ADOV-JCTB03, DDOSRSV-JCTB04, DOSV-Comb98, DOSV-JCTB00} for other results related to tree-width and vertex partitions. 

Tree-partition-width is not bounded above by any function solely of tree-width. For example, wheel graphs have bounded tree-width and unbounded tree-partition-width, as observed by \citet{BodEng-JAlg97}. Thus, it seems unavoidable that the maximum degree appears in an upper bound on the tree-partition-width. This fact, along with other applications, motivated Dujmovi\'c\etal\citep{DMW-GraphLayout,DujWoo-WG03} to study the structure of the bags in a tree-partition. In this paper we continue this approach, and prove the following result (in \secref{Main}). 

\begin{theorem}
\thmlabel{Main}
Let $k$ and $\ell$ be integers with $k\geq\ell\geq0$. Let $t=\floor{k/(\ell+1)}$. Every $k$-tree $G$ has an $\ell$-tree-partition in which each bag induces a connected $t$-tree in $G$.
\end{theorem}

It is easily seen that \thmref{Main} is tight for $G=K_{k+1}$ and for all $\ell$. Note that \thmref{Main} can be interpreted as a statement about chromomorphisms (see \citep{Strausz-PhD,Strausz03}).

Dujmovi\'c\etal\citep{DMW-GraphLayout,DujWoo-WG03} proved that every $k$-tree has a tree-partition in which each bag induces a $(k-1)$-tree. Thus \thmref{Main} with $\ell=1$ improves this result. That said, the tree-partition of Dujmovi\'c\etal\citep{DMW-GraphLayout,DujWoo-WG03} has a number of additional properties that were important for the intended application. We generalise these additional properties in \secref{Oriented}. The price paid is that each bag may now induce a $(k-\ell)$-tree, thus matching the result of Dujmovi\'c\etal\citep{DMW-GraphLayout,DujWoo-WG03} for $\ell=1$. Note that the proof of Dujmovi\'c\etal\citep{DMW-GraphLayout,DujWoo-WG03} uses a different construction to the one given here.

\section{Proof of \thmref{Main}}
\seclabel{Main}

We proceed by induction on $|V(G)|$. If $V(G)=\emptyset$, then the result holds with $V(H)=\emptyset$ regardless of $k$ and $\ell$. Now suppose that $|V(G)|\geq1$. Thus $G$ has a vertex $v$ such that $G\setminus v$ is a $k$-tree, and $N(v)$ is a $k'$-clique for some $k'\leq k$. By induction, $G\setminus v$ has an $\ell$-tree-partition $H$ in which each bag induces a connected $t$-tree. Let $C=\{x\in V(H):N(v)\cap H_{x}\ne\emptyset\}$. Since $N(v)$ is a clique, $C$ is a clique of $H$ (by the definition of $H$-partition). Since $H$ is an $\ell$-tree, $|C|\leq\ell+1$.

\textbf{Case 1.} $|C|\leq\ell$: Add one new node $y$ to $H$ adjacent to each node $x\in C$. Since $C$ is a clique of $H$ and $|C|\leq\ell$, $H$ remains an $\ell$-tree. Let $H_y=\{v\}$. The other bags remain unchanged. Since $t\geq0$, $H_y$ induces a connected $t$-tree ($=K_1$) in $G$. Thus $H$ is now a partition of $G$ in which each bag induces a connected $t$-tree in $G$. 

\textbf{Case 2.} $|C|=\ell+1$: There is a node $y\in C$ such that $|N(v)\cap H_y|\leq t$, as otherwise $|N(v)|\geq(t+1)|C|=(\floor{k/(\ell+1)}+1)(\ell+1)\geq k+1$. Add $v$ to the bag $H_y$. Let $u\in N(v)\cap H_y$. Every neighbour of $v$ not in $H_y$ is adjacent to $u$ (in $G\setminus v$). Thus $H$ is a partition of $G$. $H_y$ induces a connected $t$-tree in $G$, since $H_y\setminus\{v\}$ induces a connected $t$-tree in $G\setminus v$, and the neighbourhood of $v$ in $H_y$ is a clique of at least one and at most $t$ vertices. The other bags do not change. Thus each bag of $H$ induces a connected $t$-tree in $G$. \qed


\section{Oriented Partitions}
\seclabel{Oriented}

Let $G$ be an oriented graph with arc set $A(G)$. Let $\widehat{G}$ be the underlying undirected graph of $G$. The in- and out-neighbourhoods of a vertex $v$ of $G$ are respectively denoted by $N^-(v)=\{u\in V(G):uv\in A(G)\}$ and $N^+(v)=\{w\in V(G):vw\in A(G)\}$. It is easily seen that an (undirected) graph $G$ is a $k$-tree if and only if there is an acyclic orientation of $G$ such that for every vertex $v$ of $G$, $N^-(v)$ is a $k'$-clique for some $k'\leq k$. An oriented graph with this property is called an \emph{oriented $k$-tree}. Let $G$ and $H$ be oriented graphs. An \emph{oriented $H$-partition} of $G$ is an $\widehat{H}$-partition of $\widehat{G}$ such that for every arc $xy$ of $H$, and for every edge $vw$ of $\widehat{G}$ with $v\in H_x$ and $w\in H_y$, $vw$ is oriented from $v$ to $w$. This concept is similar to an oriented homomorphism (see \citep{Sopena-JGT97,BFKRS-JCTB01} for example).

\begin{theorem}
\thmlabel{Oriented}
Let $k$ and $\ell$ be integers with $k\geq\ell\geq0$. Let $t=k-\ell$. Every oriented $k$-tree $G$ has an oriented $\ell$-tree partition $H$ in which each bag induces a weakly connected oriented $t$-tree in $G$. Moreover, for every node $x$ of $H$, the set of vertices $Q(x)=\bigcup_{v\in H_x}(N^-(v)\setminus H_x)$ is a $k'$-clique of $G$ for some $k'\leq k$. 
\end{theorem}

The construction in the proof of \thmref{Oriented} only differs from that of \thmref{Main} in the choice of the node $y$ in Case~2.

\begin{proof}
We proceed by induction on $|V(G)|$. If $V(G)=\emptyset$, then the result holds with $V(H)=\emptyset$ regardless of $k$ and $\ell$. Now suppose that $|V(G)|\geq1$. Since $G$ is acyclic, there is a vertex $v$ of $G$ such that $N^+(v)=\emptyset$, $N^-(v)$ is a $k'$-clique for some $k'\leq k$, and $G\setminus v$ is an oriented $k$-tree. By induction, there is an oriented $\ell$-tree-partition $H$ of $G\setminus v$ in which each bag induces a weakly connected oriented $t$-tree in $G\setminus v$. Moreover, for every node $x$ of $H$, $Q(x)$ is a $k'$-clique for some $k'\leq k$. Let $C=\{x\in V(H):N^-(v)\cap H_x\ne\emptyset\}$. Since $N^-(v)$ is a clique, $C$ is a clique of $H$. Since $H$ is an oriented $\ell$-tree, $|C|\leq\ell+1$.

\textbf{Case 1.} $|C|\leq\ell$: Add one new node $y$ to $H$ adjacent to each node $x\in C$. Orient each new edge from $x$ to $y$. Obviously $H$ remains acyclic. Since $C$ is a clique of $H$ and $|C|\leq\ell$, $H$ remains an oriented $\ell$-tree. Let $H_y=\{v\}$. The other bags are unchanged. Since $t\geq0$, $H_y$ induces a weakly connected oriented $t$-tree ($=K_1$) in $G$. All edges of $G$ that are incident to a vertex in $H_y$ are oriented into the vertex in $H_y$. Thus $H$ is now an oriented partition of $G$ in which each bag induces a weakly connected oriented $t$-tree in $G$. Now $Q(y)=N^-(v)$, which is a $k'$-clique for some $k'\leq k$. $Q(x)$ is unchanged for nodes $x\ne y$. Hence the theorem is satisfied. 

\textbf{Case 2.} $|C|=\ell+1$: The clique $C$ induces an acyclic tournament in $H$. Let $y$ be the sink of this tournament. Since $|N^-(v)\cap H_x|\geq1$ for every node $x\in C\setminus\{y\}$, $|N^-(v)\cap H_y|\leq k'-(|C|-1)\leq k-\ell=t$. Add $v$ to the bag $H_y$. 

Consider a neighbour $u$ of $v$. Since $N^+(v)=\emptyset$, $uv$ is oriented from $u$ to $v$. Say $u\in H_z$ with $z\ne y$. Then $z$ is in the clique $C$. Thus $zy$ is an edge of $H$. Since $y$ is a sink of $C$, $zy$ is oriented from $z$ to $y$. Thus $H$ is now an oriented partition of $G$. $H_y$ induces a weakly connected oriented $t$-tree in $G$, since $H_y\setminus\{v\}$ induces an oriented $t$-tree in $G\setminus v$, and the in-neighbourhood of $v$ in $H_y$ is a clique of at least one and at most $t$ vertices. The other bags do not change. Thus each bag of $H$ induces a weakly connected oriented $t$-tree in $G$. 

$Q(y)$ is not changed by the addition of $v$ to $H_y$, as there is at least one vertex $u\in N^-(v)\cap H_y$, and any vertex in $N^-(v)\setminus H_y$ is also in $N^-(u)\setminus H_y$. For nodes $x\ne y$, $Q(x)$ is unchanged by the addition of $v$ to $H_y$, since $v$ is not in the in-neighbourhood of any vertex. Hence the theorem is satisfied. 
\end{proof}

\section*{Acknowledgements} 

Thanks to Matthew DeVos, Vida Dujmovi\'c, Attila P\'or, and Ricardo Strausz  for stimulating discussions.


\def\cprime{$'$} \def\soft#1{\leavevmode\setbox0=\hbox{h}\dimen7=\ht0\advance
  \dimen7 by-1ex\relax\if t#1\relax\rlap{\raise.6\dimen7
  \hbox{\kern.3ex\char'47}}#1\relax\else\if T#1\relax
  \rlap{\raise.5\dimen7\hbox{\kern1.3ex\char'47}}#1\relax \else\if
  d#1\relax\rlap{\raise.5\dimen7\hbox{\kern.9ex \char'47}}#1\relax\else\if
  D#1\relax\rlap{\raise.5\dimen7 \hbox{\kern1.4ex\char'47}}#1\relax\else\if
  l#1\relax \rlap{\raise.5\dimen7\hbox{\kern.4ex\char'47}}#1\relax \else\if
  L#1\relax\rlap{\raise.5\dimen7\hbox{\kern.7ex
  \char'47}}#1\relax\else\message{accent \string\soft \space #1 not
  defined!}#1\relax\fi\fi\fi\fi\fi\fi}

\end{document}